\documentclass[12pt]{article}
\usepackage[latin1]{inputenc}
\usepackage{amsfonts}
\usepackage{indentfirst}
\usepackage{amssymb,amsmath}
\usepackage{graphicx}
\newtheorem{Theorem}{Theorem} 
\newtheorem{Definition}{Definition} 
\newtheorem{Proposition}{Proposition} 
\newtheorem{Lemma}{Lemma} 
\newtheorem{Corollary}{Corollary} 
\newtheorem{Exam}{Example} 

\newcommand{\bdem}{\textsc{Proof}\rm}
\newcommand{\edem}{\hfill $\Box$ \\[3mm]}
\title{A partial reciprocal of Dirichlet Lagrange Theorem detected by Jets}
\author{Alva, G.~J.~M. and Garcia, M. V.~P.\footnote{The first author was partially supported by CNPq-Brasil. The second author was partially supported by FAPESP.}}

\date{}

\begin{document}
\maketitle
\begin{abstract}
We study the stability of an equilibrium point in a conservative Hamiltonian system in the case that equilibrium is not a minimum of the potential energy and this fact is shown by a jet of this function. Thanks to a modification of a result of Krasovskii, we prove that for a large class of systems under these conditions equilibrium is unstable and there is an asymptotic trajectory to that point.
\end{abstract}

\section*{Introduction}\label{sec1}
We study in this work the stability of an equilibrium point of a Hamiltonian system of $n$ degrees of freedom, sufficient smooth, of
potential energy $\Pi:\Omega\rightarrow \mathbb{R}$, where $\Omega$ is an open neighborhood of the origin of $\mathbb{R}^{n}$, and
cinetic energy $T:\Omega\times\mathbb{R}^{n} \rightarrow \mathbb{R}$, $T(q,p)=\frac{1}{2}\langle B(q)p,p\rangle$, in that for all 
$q\in \Omega$, $B(q)$ is a matrix $n\times n$, symmetric positive definite. In these conditions, if the origin is a critical point of
$\Pi$, $(0,0)$ is a equilibrium point of the hamiltonian field defined in $\Omega\times\mathbb{R}^{n}$ of equations
\begin{eqnarray}\label{eq1sec1}
X_{H}(q,p)\left \{\begin{array}{lll}
\dot{q}=\frac{\partial H}{\partial p}(q,p)\\
&&\\
\dot{p}=-\frac{\partial H}{\partial q}(q,p)
\end{array}\right.
\end{eqnarray}

As is well known, the Dirichlet-Lagrange theorem guarantees that, if the origin is a local strict minimum point of the potential energy
then $(0,0)$ is a stable equilibrium in the sense of Lyapunov of (\ref{eq1sec1}), this result does not admit in the strict sense of the
term, a reciprocal, because as shown by Painleve (see \cite{PP}) in the early twentieth century, if $\Pi:\mathbb{R}\rightarrow \mathbb{R}$ is
$\Pi(q)=e^{-1/q^{2}}\sin(1/q)$, if $q\neq 0$, $\Pi(0)=0$ e $T(q,p)=\frac{p^{2}}{2}$, $(0,0)$ is a stable equilibrium of (\ref{eq1sec1}),
although $0$ is not a point of minimum of $\Pi$.

Since then the question of determining sufficient conditions on the potential energy to ensure instability of $(0,0)$ when this is not a point of minimum of $\Pi$ is known in the literature as ``the problem of inversion of the Dirichlet-Lagrange theorem''
 and has attracted the attention
of mathematicians as Lyapunov, Cetaev, La Salle, Salvadori, Kozlov, Laloy, among others (see...for an historical overview of this problem)

Our work is connected to two basic ideas, Lyapunov proof in \cite{L} that if the Hessian of $\Pi$ at the origin shows that this is not a
point of minimum of the potential energy then $(0,0)$ is an unstable equilibrium de (\ref{eq1sec1}), the Lyapunov own launched the idea
of generalizing this result, conjecturing that if the Taylor polynomial of order $s$ in $0$ guarantee that this is not point of minimum
of the potential energy then $(0,0)$ would be unstable.

Using the language of punctual jets, introduced by A. Barone in \cite{B}, this Lyapunov conjecture says that \emph{if the punctual jet of order
$s$ of $\Pi$ at $0$ shows that the potential energy $\Pi$ has no minimum in this point, then the equilibrium $(0,0)$ of
(\ref{eq1sec1}) is unstable}.

This conjecture was investigated by Maffei, Moauro and Negrine in two notable works, \cite{MN} and \cite{MMN}, they study the case where the jet showing that the origin is not a minimum, is
an homogeneous polynomial. In fact, they prove more, the main result of these articles proves instability of the origin when, after an eventual change of coordinates, the potential energy has the form $\pi(q_{1},\ldots,q_{n})=\Vert{(q_{1},\ldots,q_{r})}\Vert^{2}+\pi_{s}(q_{r+1},\ldots,q_{s})+R(q_{r+1},\ldots,q_{s})$, where $0\le r<s$, $\pi_{s}$ is an homogeneous polynomial of degree $s\ge 2$ which does not have a minimum at the origin, and $R$ obeys $\lim\limits_{q\to 0}\frac{R(q)}{\Vert{q}\Vert^{s}}=0$. 


Other relevant works on this conjecture were made by Freire, Garcia and Tal in \cite{GT} and \cite{FGT}, in the first of these articles the conjecture is proved in the case of systems with two degrees of freedom, the second one studies the general case n degrees of freedom and it is shown that if the punctual jet of order
$s$ of $\Pi$ at $0$,
$j^{s}\Pi$, shows that $\Pi$ have not minimum at origin and the Hessian of $\Pi$ in $0$ is
positive semi-definite and the dimension of its null space is at most $2$, we have the instability of the origin.



In the four aforementioned works the instability is proved showing that there is an asymptotic trajectory to the origin.


On the other hand, a classical result of Cetaev (see \cite{LHR}) establishes the instability of the origin when $\langle \nabla \Pi(q),q\rangle<0$, for all $q$ in a connected component of $\Pi^{-1}(-\infty,0)$. (If $\Pi$ has this property we will
say that $\Pi$ satisfies the Cetaev condition, or that $\Pi$ is a Cetaev potential).

As a result, several instability results were obtained  establishing conditions to ensure that $\Pi$ satisfies Cetaev's condition. The reader may find in section 3 of chapter 3 of the excellent book [LHR] several examples of this type of result.


Our work connects these two lines, more precisely, it is assumed that $ j^{s}\Pi $ shows that the origin is not a minimum point of
potential energy and $ j^{s}\Pi $ satisfies Cetaev's condition. These assumptions do not guarantee that $\Pi$ obeys the condition of Cetaev, so we seek to establish additional conditions on $ j^s\Pi $ to make sure this happens.


Our main results are in section \ref{sec:instbiltres}. In Theorem \ref{d5} we present a set of conditions that guarantee that if $j^{s}\Pi$ is a Ceatev potential the same occurs to $\Pi$. More precisely, there are essentially two conditions imposed on the potential energy, one of them has a topological character and is simple to see that it is natural in this context, it requires that closure of the region where $j^{s}\Pi$ is negative is contained within the region in which  $\langle\nabla j^{s}\Pi(q) | q\rangle$ is negative, the second condition requires that $j^{\ell}\Pi$ is semidefinite positive for all $\ell<s$, has a technical character that we will duscuss in final section of the article. Moreover, in the Proposition \ref{prop:karsovskii} we generalize a result of Krasoviskii about auxiliary functions and assymptotic trajectories in differential equations, and use this Proposition to show that in a great class of Cetaev potentials, which contains those who satisfy the conditions of Theorem \ref{d5}, we have an asymptotic trajectory to the equilibrium. To do it, we show in Proposition \ref{cor:hassymptotic} that, if $\langle\nabla\Pi(q) | q\rangle < 0$ for all $q\not=0$ in the closure of a connected component of $\Pi^{-1}(]-\infty,0[)$, there is an asymptotic trajectory to the equilibrium.
 
We finish this introduction with an outline of the article. In section \ref{sec2} we introduce the basic notions and terminology, along with the statement of
some useful Propositions, whose proofs we do later, in section \ref{sec.proofsomited}. The core of this work is the section \ref{sec:instbiltres} where we present Theorem \ref{d5} and Proposition \ref{cor:hassymptotic}. Finally, in section \ref{sec:discuss} we discuss the
aforementioned \emph{technical character} that we have in part of hypotheses in Theorem \ref{d5} and we give an example showing that, in order to ensure that $ \Pi $ is a Cetaev potential based only in informations about $j^{s}\Pi$, perhaps this restriction could be slightly weakened, but not much.


\section{Statement of the problem and basic results}\label{sec2}
\subsection{Preliminary notions}
We present here the definitions and notations that will be used throughout the article.
\begin{Definition}
Let $s$ be a natural nonzero and $\Omega$ an open subset of $\mathbb{R}^{n}$ containing the origin.
It is said that $f:\Omega \rightarrow \mathbb{R}$ has punctual jet of order $s$ at the origin if exists a polynomial
$p:\mathbb{R}^{n}\rightarrow \mathbb{R}$ of degree less than or equal to $s$ such that
$\lim_{x\rightarrow 0}\frac{f(x)-p(x)}{|x|^{s}}=0$. In this case we denoted $p$ by $j^{s}f$.
\end{Definition}
It is simple to see that there is at most a polynomial that satisfies these conditions, and when $j^{s}f$ exists we will simply
say that this is the $s$-jet de $f$, instead of $s$-jet de $f$ at the origin.

Is also easy to see that if $f$ is of class $\mathcal{C}^{s}$ the Taylor polynomial of order $s$ of $f$ at the origin
is the $j^{s}f$, the function $f:\mathbb{R}\rightarrow \mathbb{R}$, $f(x)=e^{-1/x^{2}}$, if $x\notin \mathbb{Q}$ and
$f(x)=0$ if $x\in \mathbb{Q}$ show that the reciprocal of this statement is not true, because for all $s\in \mathbb{N}$,
$j^{s}f\equiv 0$ e $f$ is not $\mathcal{C}^{s}$.

If $j^{s}f$ is the s-jet of $f$ and $1\le k\le s$ then $j^{k}f=j^{k}(j^{s}f)$ and, if $f_{k}=j^{k}f-j^{k-1}f$,
we see that $f_{k}$ is the null polynomial or is an homogeneous polynomial of degree $k$, we will refer to $f_{k}$ as the \emph{homogeneous part of degree $k$} of $j^{s}f$ and it's clear that $j^{s}f=\sum\limits_{k\le s}f_{k}$.

This notion was introduced by Barone-Netto in [B].\\
\begin{Definition}
Suppose that $f:\Omega \rightarrow \mathbb{R}$ has $s$-jet. It is said that $j^{s}f$ shows that $f$ has not minimum at
the origin if for every $g:\Omega \rightarrow \mathbb{R}$ with $j^{s}g=j^{s}f$ we have that $g$ has no minimum (even a weak minimum) at the origin.
\end{Definition}
We work with functions $f:\Omega\stackrel{\mathcal{C}^{2}}{\rightarrow}\mathbb{R}$ that have $s$-jet and we will call
the function $R_{f}(q)=\langle \nabla f(q),q\rangle,\ q\in\Omega$ as the \emph{radial derivative} of $f$. 

If $f$ has $s$-jet, the radial derivative of $j^{s}f$ will be denoted by $R_{f}^{s}$.
When it is clear what is the function $f$ involved we will use $R(q)$ and $R^{s}(q)$ instead of $R_{f}(q)$ and $R_{f}^{s}(q)$.

\subsection{The problem}
Consider an Hamiltonian system of potential energy $\Pi:\Omega\stackrel{\mathcal{C}^{2}}{\rightarrow}\mathbb{R}$, where $\Omega\subset\mathbb{R}^{n}$ is an open neighbourhood of the origin and kinetic energy $T:\Omega\times\mathbb{R}^{n}\stackrel{\mathcal{C}^{2}}{\rightarrow}\mathbb{R}$, $T(q,p)=\frac{1}{2}<B(q)p|p>$, where $B(q)$ is a symmetric positive definite matrix.   

Suppose that origin is a critical point of the potential energy, $\Pi(0)=0$ and admit that $\Pi$ has $s$-jet at origin, so  $\Pi(q)=j^{s}\Pi(q)+o(|q|^{s})$. Moreover, $j^{1}\Pi$ is null, then  $j^{s}\Pi=\sum\limits_{\ell=2}^{s}\Pi_{\ell}$.

Remember that $R^{s}_{\Pi}(q)=\langle \nabla j^{s}\Pi(q),q\rangle$, consider $A_{s}= (j^{s}\Pi)^{-1}((-\infty,0))$, $C_{s}= (R_{\Pi}^{s})^{-1}((-\infty,0))$ and denote the open ball of center $0$ and radius $r$ by $B_{r}$. Also, denote the interior of a subset $U$ of $\mathbb{R}^{n}$ by $U^{\circ}$.

In this conditions $(0,0)$ is an equilibrium of the equations of the system, namely equations \eqref{eq1sec1}. 

We suppose that $j^{s}\pi$ has the Cetaev's condition, i.e. $R^{s}_{\Pi}(q)<0$ if $j^{s}\Pi(q)<0$, and we will investigate under what conditions this ensures that $\Pi$ also satisfies Cetaev's condition.

In order to do it, we admit that potential energy obeys
the following conditions:
\begin{enumerate}
\item [$(H_{1})$] $j^{\ell}\Pi(q)\geq 0$ in some neighborhood of the origin if $\ell\leq s-1$;
\item [$(H_{2})$] $j^{s}\Pi$ show that $\Pi$ has not minimum at origin;
\item [$(H_{3})$] There exist $\varepsilon > 0$ with $\overline{B_{\varepsilon}}\subset \Omega$ and a connected component $A_{*}$ of $A_{s}$, such that
$(\overline{A}_{*}\setminus \{0\})\cap B_{\varepsilon}\subset C_{s}^{\circ}\cap B_{\varepsilon}$.
\end{enumerate}

It is said that $\Pi$ is of class $\mathcal{H}$ if $\Pi$ satisfies the hypotheses $(H_{1})$-$(H_{3})$.

Note that by stating that $j^{s}\Pi$ satisfies the condition of
Ceatev, one is saying that $A_{s}\subset C_{s}$, while $(H_{3})$
establishes a stronger condition. However, the condition $A_{s}\subset
C_{s}$ does not guarantee that $\Pi$ obeys the Cetaev's condition,
even when $(H_{1})$ and $(H_{2})$ are valid. To see this, just look at
the example $\pi(x,y)=y^{2}-x^{2}+x^{3}$, it's clear that $j^{2}\pi$
satisfies Cetaev's condition and obeys $(H_{1})-(H_{2})$, but $\pi$
has not Cetaev's condition. This shows that condition $(H_{3})$ is a
natural one in this context.  

We will say that a potential energy $\Pi$ has the \emph{strict Cetaev's condition} if there is a connected component $A$ of $\Pi^{-1}(]-\infty,0[)$ such that $R_{\Pi}(q)<0$, for $q\in\overline{A}\setminus\{0\}\cap B_{\varepsilon}$, in other words, if $\overline{A}\setminus\{0\}\cap B_{\varepsilon}\subset R_{\Pi}^{-1}(]-\infty,0[)$.

In the main results of this paper we will prove that if the potential
energy is of class $\mathcal{H}$ then it satisfies Cetaev's strict 
condition (Theorem \ref{d5}), and we show that if $\pi$ has the Cetaev's strict 
condition, then 
there is an asymptotic trajectory to the origin in equations
\eqref{eq1sec1} (see Proposition \ref{cor:hassymptotic}).   


\subsection{Basic results}\label{sec.statements}
Here we state the technical results needed for demonstration of our theorems of instability, presented in the next section. The proofs of the propositions set out in this section without references are found in section \ref{sec.proofsomited}.

We start by the observation that an algebraic curve  $\gamma :[0,\rho]\rightarrow \mathbb{R}^{n}$, with  $\gamma(0)=0\in \Omega$ and $\gamma(t)\neq 0$ for all
$t\in (0,\rho]$, is tangent at origin to a half-line. We will denote this half-line by $r_{\gamma}$ (a reference to see this and other tecnical facts about algebraic curves is \cite{JM}).

A very usefull consequence of hypothesis $(H_{3})$ is about the geometry of the set $A_{s}$ and the tangent half-lines of algebraic curves in $A_{s}$ that starts at origin, precisely:    
\begin{Lemma}\label{b1}
	Suppose that $\Pi$ satisfies $(H_{3})$ and let be $\gamma :[0,\rho]\rightarrow \mathbb{R}^{n}$
	an algebraic curve, with $\gamma(0)=0$ and $\gamma(t)\in \overline{A}_{s}\setminus \{0\}$ if $t>0$. Then
	$j^{s}\Pi|_{r_{\gamma}}$ has a strict local maximum at the origin, furthermore there is $\varepsilon>0$ such that $r_{\gamma}\setminus \{0\}\cap B_{\varepsilon}$ and the image of $\gamma$ are contained in the same connected component of $\overline{A}_{s}$.
\end{Lemma}
In terms of connected components this result has a very usefull form.
\begin{Corollary}\label{cor:compconexa}
	Suppose that there are connected components $A_{*}$ of $A_{s}$ and $C_{*}$ of $C_{s}$ such that, for some $\varepsilon>0$, $(\overline{A_{*}}\setminus\{0\})\cap B_{\varepsilon}\subset C_{*}\cap B_{\varepsilon}$.
	Then, if  $\gamma :[0,\rho]\rightarrow \mathbb{R}^{n}$ is
	an algebraic curve, with $\gamma(0)=0$ and $\gamma(t)\in \overline{A}_{*}\setminus \{0\}$ if $t>0$, we have that $r_{\gamma}\setminus\{0\}\cap B_{\varepsilon}$ is contained in $A_{*}$.
\end{Corollary}
In view of the hypotheses $\mathcal{H}$ the set $Z_{s-1}=\cap_{\ell=2}^{s-1}\{q\in \Omega : \Pi_{\ell}(q)= 0\}\setminus \{0\}$ will have a great importance to us. Since $\Pi_{\ell}$ are homogeneous polynomials, if $q\in Z_{s-1}$ and $\lambda>0$ then $\lambda q\in Z_{s-1}$. We will call a subset of $\mathbb{R}^{n}$ with this property of cone of vertex at origin.

We will now construct a cone with vertex 0, vitally important for our
purposes. As a matter of convenience in the reading, we present the
properties of this cone in three stages, covered the following three
propositions.
\begin{Proposition}\label{c2}
Suppose that $\Pi$ satisfies the hypothesis $\mathcal{H}$, then there
is a half-line $r$ of origin $0$ contained in $Z_{s-1}$ such that
$\Pi_{s}(q)<0$ for all $q\in r$. Moreover, if $\Delta_{r}$ is the
connected component of $Z_{s-1}$ which contains $r$, there exists a
cone 
$K_{s}\subset \mathbb{R}^{n}$ with vertex $0$ such that 
\begin{enumerate}
\item[$(a)$] $(\Delta_{r}\setminus\{0\})\subset (K_{s})^{\circ}$ 
\item[$(b)$] If $q\in \partial K_{s}\setminus \{0\}$ there exists $\ell \in \{2,\cdots, s-1\}$ such that $\Pi_{\ell}(q)>0$ and
$\Pi_{j}(q)=0$ for $j< \ell$.
\item[$(c)$] $\Pi_{s}(q)<0$ for all $q\in K_{s}\setminus \{0\}$
\end{enumerate}
\end{Proposition}

\begin{Proposition}\label{c3}
Suppose that $\Pi$ satisfies the hypothesis $\mathcal{H}$
and consider the cone $K_ {s}$ given in the previous proposition, then
there are constants $\varepsilon>0$ and $m>0$ such that 
\begin{enumerate}
\item[$(1)$]
$j^{s-1}\Pi(q)\geq \frac{m}{2}|q|^{s-1}$, for each $q\in (\partial K_{s}\setminus\{0\})\cap B_{\varepsilon}$
\item [$(2)$]
$R_{\Pi}^{s-1}(q)\geq (s-1)\frac{m}{2}|q|^{s-1}$, for each $q\in (\partial K_{s}\setminus\{0\})\cap B_{\varepsilon}$
\end{enumerate}
\end{Proposition}
At last we have
\begin{Proposition}\label{c5}
Suppose that $\Pi$ satisfies the hypothesis $\mathcal{H}$, then there is $\varepsilon>0$ such that the cone $K_{s}$ given in the above propositions obeys  
\begin{enumerate}
\item[$(i)$] $j^{s}\Pi(q)>0$ for each $q\in (\partial K_{s}\setminus\{0\})\cap B_{\varepsilon}$
\item[$(ii)$] $R_{\Pi}^{s}(q)>0$ for each $q\in (\partial K_{s}\setminus\{0\})\cap B_{\varepsilon}$
\end{enumerate}
\end{Proposition}
As a consequence of proposition \ref{c2} is clear that $M_{s-1}:=\left\{q\in Z_{s-1}: \Pi_{s}(q)<0\right\}$ is a nonempty set. 

 The set $M_{s-1}$ is indeed a
 positive cone with vertex at the origin and, for any
 algebraic curve $\gamma :[0,\rho]\longrightarrow \Omega $ with
 $\gamma(0)=0$ and $\gamma(t)\in \overline{A}_{s}\setminus \{0\}$ if
 $0< t \leq \rho$, it follows from lemma \ref{b1} that
 $r_{\gamma}\subset M_{s-1}$. Because of this $M_{s-1}$ will be called
 the \emph{tangent cone} of $\overline{A}_{s}$.



\begin{Proposition}\label{c7}
Suppose that $\Pi$ satisfies the hypothesis $\mathcal{H}$, and let $r$
be a half-line with origin $0$ contained in $M_{s-1}$. If $A_{*}$ and
$C_{*}$ are, respectively the connected component of $A_{s}$ and
$C_{s}$ containing $r\setminus \{0\}$ and $K_{s}$ is a cone given by
the proposition \ref{c2}, then there are $\varepsilon_{0}>0$ and
$M_{2}>0$ such that
\begin{eqnarray*}
	R_{\Pi}^{s}(q) \le -M_{2}|q|^{s},\quad \forall q\in (\partial A_{*}\setminus \{0\})\cap B_{\varepsilon_{0}}
\end{eqnarray*}
and
\begin{eqnarray*}
(\overline{C}_{*}\setminus \{0\})\cap B_{\varepsilon_{0}}\subset (K_{s})^{\circ}\cap B_{\varepsilon_{0}}. 
\end{eqnarray*}
\end{Proposition}
For convenience, from now on we will call a cone $K_ {s}$ that satisfies the propositions 1-4 simply cone built in the section 2

\section{Instability results}\label{sec:instbiltres}

Here we show that if the potential energy satisfies the hypotheses
$\mathcal{H}$  then $\Pi$ is a Cetaev's potential. So Cetaev's theorem
immediately implies that $\mathcal{H}$ causes the instability of
$(0.0)$ in equations (\ref{eq1sec1}). Actually we will prove later
that, under these conditions, there is an asymptotic trajectory to the
origin. 

We begin by noting that, by Euler's theorem on homogeneous functions,
one can see that 

\begin{equation}\label{eq:euler}
R_{\Pi}^{s}(q)=(s-1)j^{s}\Pi(q)-\sum_{\ell=2}^{s-2}j^{\ell}\Pi(q)+\Pi_{s}(q).
\end{equation}


\begin{Proposition}\label{d1}
	Suppose that $\Pi$ satisfies the hypotheses $\mathcal{H}$, let
        $r$ be a half-line with origin $0$ contained in $M_{s-1}$, and
        denote $A_{*}$ and $C_{*}$, respectively, the  connected
        components of $A_{s}$ and  $C_{s}$ containing $r\setminus \{0\}$.
	Then there exists $0<\varepsilon<1$ and a semi-algebraic set
        $W_{s}\subset B_{\varepsilon}$, such that $0\in \partial
        W_{s}$, and:  
\begin{enumerate}
\item $(\overline{A}_{*}\setminus \{0\})\cap B_{\varepsilon} \subset W_{s}^{\circ}\subset \overline{W}_{s}\setminus \{0\} \subset (C_{*})^{\circ}\cap B_{\varepsilon}$
\item $\Pi_{s}(q)<0$ for all $q\in \overline{W}_{s}\setminus \{0\}$ and, moreover, there are constants $\alpha>0$ and $\beta>0$ such that:
\begin{enumerate}
\item[$(a)$]If $q\in \partial W\setminus \{0\}$, then $j^{s}\Pi(q)\geq -\alpha \Pi_{s}(q)$;
\item[$(b)$]If $q\in \overline{W}\setminus \{0\}$, then $R_{\Pi}^{s}(q)\leq \beta \Pi_{s}(q)$
\end{enumerate}
\end{enumerate}
\end{Proposition}
\bdem\ Consider a cone $K_{s}$ built in section 2 and
$0<\varepsilon<\varepsilon_{0}<1$, whith $\varepsilon_{0}$ given by
proposition \ref{c7}.  

Let $Q$ be the polynomial
$Q(q)=(s-1)j^{s}\Pi(q)+\frac{1}{2}\left[-\sum_{\ell=2}^{s-2}j^{\ell}\Pi(q)+\Pi_{s}(q)\right]$
and consider the semi-algebraic set $W=\{q\in
B_{\varepsilon}:Q(q)<0\}$.

Note that $0\in\partial W$, since $Q(q)<0$ for $q\in r$.

It's clear from $\mathcal{H}$ that $\overline{A}_{*}\setminus
\{0\}\cap B_{\varepsilon} \subset W $. 

Moreover, by proposition \ref{c2}, $\Pi_{s}(q)<0$, if $q\in
K_{s}\setminus\{0\}$ , then by using \eqref{eq:euler} one can see that
$R_{\Pi}^{s}(q)=Q(q)+\frac{1}{2}\left[-\sum_{\ell=2}^{s-2}j^{\ell}\Pi(q)+\Pi_{s}(q)\right]$
and conclude $W\cap K_{s} \subset C_{*}\cap B_{\varepsilon}$. Since,
by proposition \ref{c7}, $\overline{C}_{*}\setminus \{0\}\cap
B_{\varepsilon_{0}}\subset (K_{s})^{\circ}$, it follows from here that
$\overline{W}\setminus\{0\}\subset C_{*}\cap B_{\varepsilon}\subset
(K_{s})^{\circ}$, which concludes the proof of $1$. 

So $\Pi_{s}(q)<0$ for all $q\in\overline{W}\setminus\{0\}$, then,
if $q\in\partial W\setminus\{0\}$, we have 
$Q(q)=0$, i.e.
$(s-1)j^{s}\Pi(q)=-(1/2)\left[-\sum_{\ell=k}^{s-2}j^{\ell}\Pi(q)+\Pi_{s}(q)\right]\geq 
-(1/2)\Pi_{s}(q)$, which implies 2-a for $\alpha=\frac{1}{2(s-1)}$.  

On the other hand, if $q\in \overline{W}\setminus \{0\}$, then
$Q^{s}(q)\leq 0$, getting the second inequality with
$\beta=\frac{1}{2}$ since, by \eqref{eq:euler}, $R_{\Pi}^{s}(q)\leq
(1/2)\left[-\sum_{\ell=k}^{s-2}j^{\ell}\Pi(q)+\Pi_{s}(q)\right]\leq
(1/2)\Pi_{s}(q)$. 
\edem



Now we use the homogeneity of $\Pi_{s}$ and the above proposition to
get the crucial result of this section.

\begin{Proposition}\label{d3}
Suppose $ \Pi $ satisfies the hypotheses $ \mathcal {H} $, and let $W$ be
the set obtained in proposition \ref{d1}, then there is
$0<\varepsilon_{1}<1$ such that: 
\begin{enumerate}
\item[$(i)$]
$\Pi(q)>0$ for all $q\in
(\partial W\setminus \{0\})\cap B_{\varepsilon_{1}}$. 
\item[$(ii)$]
$R_{\Pi}(q)<0$ for all
$q\in\overline{W}\setminus \{0\})\cap B_{\varepsilon_{1}}$.
\end{enumerate}
\end{Proposition}
\bdem\ Let $g(q)=\Pi(q)-j^{s}\Pi(q)$ and
$h(q)=R_{\Pi}(q)-R^{s}_{\Pi}(q)$, clearly $g$ and $h$ are functions
$o(\vert{q}\vert^{s})$ at origin.

Consider the cone $K_{s}$ built in section 2 and take the
semi-algebraic set $W$ and the constants
$0<\varepsilon<1,\ \alpha>0, \beta>0$ obtained in proposition
\ref{d1}. 

Then
$\overline{W}\setminus \{0\}\subset(K_{s})^{\circ}\cap
B_{\varepsilon}$.

So, by proposition \ref{c2} $\Pi_{s}(q)<0$, for all
$q\in\overline{W}\cap S_{\varepsilon}$, which is a compact set, so
$\min \{\Pi_{s}(q):q\in\overline{W}\cap S_{\varepsilon}\}=-m<0$.

By the homogeneity of $\Pi_{s}$ and
$\overline{W}\subset(K_{s})^{\circ}$, it follows that
$\Pi_{s}(q)\le-m\vert{q}\vert^{s}$, for all $q\in W$.

Using that $g$ and $h$ are $o(\vert{q}\vert^{s})$ at origin, it is
immediate to see that there is $0<\varepsilon_{1}\le\varepsilon$ such
that $\vert{g(q)}\vert\le\frac{\alpha m}{2}\vert{q}\vert^{s}$  and  
$\vert{h(q)}\vert\le\frac{\beta m}{2}\vert{q}\vert^{s}$ for all $q\in
B_{\varepsilon_{1}}$. 

Then, by proposition \ref{d1}, if
$q\in\partial W\setminus\{0\}\cap B_{\varepsilon_{1}}$, we have
\begin{center}
  $\frac{\Pi(q)}{\vert{q}\vert^{s}}=
  \frac{j^{s}\Pi(q)+g(q)}{\vert{q}\vert^{s}}\ge
  \frac{-\alpha\Pi_{s}(q)+g(q)}{\vert{q}\vert^{s}} \ge 
  \frac{\alpha m}{2}>0$
\end{center}
and, if $q\in W\setminus\{0\}\cap B_{\varepsilon_{1}}$,
\begin{center}
  $\frac{R_{\Pi}(q)}{\vert{q}\vert^{s}}=
  \frac{R_{\Pi}^{s}(q)+h(q)}{\vert{q}\vert^{s}}\le
  \frac{\beta\Pi_{s}(q)+h(q)}{\vert{q}\vert^{s}} \le \frac{-\beta
    m}{2}<0$,
\end{center}
wich ends the proof. \edem

With this proposition, it is easy to prove that if  $\Pi$ satisfies $\mathcal{H}$
then it is a Cetaev's potential, as we claim before. 


\begin{Theorem}\label{d5}
If the potential energy of a conservative Hamiltonian system is of class $\mathcal{H}$ then it has the strict Cetaev's condition.
\end{Theorem}
\bdem\ As before, consider $Z_{s-1}=\cap_{\ell=2}^{s-1}\{q\in\Omega:\Pi_{\ell}(q)=0\}\setminus\{0\}$, let $r$ be a ray of origin $0$ contained in $Z_{s-1}$ with $\Pi_{s}(q)<0$, for all $q\in r$. 

Then there is $\varepsilon>0$ such that $\Pi(q)<0$ and $R_{\Pi}(q)<0$, for all $q\in r\cap B_{\varepsilon}$, it's clear that we can suppose $0<\varepsilon<\varepsilon_{1}$, where $\varepsilon_{1}$ is given in proposition \ref{d3}.

Now consider $W$ obtained in proposition \ref{d1} and denote $A$ and $C$ respectively, the connected components of $\Pi^{-1}(]-\infty,0[)$ and $R_{\Pi}^{-1}(]-\infty,0[)$ that contains $r\cap B_{\varepsilon}$.

By proposition \ref{d3} it follows that $\pi(q)>0$, for all $q\in\partial W\setminus\{0\}\cap B_{\varepsilon}$, and $R_{\Pi}(q)<0$ for all $q\in\overline{W}\setminus\{0\}\cap B_{\varepsilon}$.

Then $A\cap B_{\varepsilon}\subset\overline{W}\cap B_{\varepsilon}\subset C\cap B_{\varepsilon}$. \edem


\subsection{Asymptotic trajectories}\label{sec.assinptotic}
The instability of origin for Hamiltonian systems whose potential energy is of class $\mathcal{H}$ follows immediately from last result, but as we have said, one can prove that in these conditions there is an asymptotic trajectory, and it's what we will do now, using a theorem inspired in a Krasovskii result. 

Consider the differential equation $\dot x=f(x)$, where $\Delta\subset\mathbb{R}^{n}$ is an open neighbourhood of the origin, $f:\Delta\longrightarrow\mathbb{R}^{n}$ is a function of class $\mathcal{C}^{1}$, with $f(0)=0$. 

The aforementioned Krasovskii result (see theorem 12.1 in \cite{K}) is a particular case of next statement when $V=W$.
\begin{Proposition}\label{prop:karsovskii}
	Suppose that $U\subset\mathbb{R}^{n}$ is an open set with $0\in\partial U$, and let $V,\ W$, $\mathcal{C}^{1}$ real functions defined in $\Delta$, such that:
	\begin{enumerate}
		\item[(i)] $\dot{V}(x)<0$, if $x\in\overline{U}\setminus\{0\}$.
		\item[(ii)] $W(x)<0$, if $x\in U$,
		\item[(iii)] $W(x)=0$, if $x\in \partial U$,
		\item[$(iv)$] $\dot{W}(x)\leq 0$, if $ x\in U$.
	\end{enumerate}
	Then there exists a trajectory $x(t)$ of $\dot{x}=f(x)$ such that $\lim\limits_{t\to\-\infty}x(t)=0$.
\end{Proposition}
\bdem\ 
Take a $\varepsilon>0$ such that $\overline{B}_{\varepsilon}\subset\Delta$
and denote by $U_{\varepsilon}$ the set $B_{\varepsilon}\cap \overline{U}$.

Consider a point $z\in U_{\varepsilon}^{\circ}$.

We'll prove first that the solution $x_{z}$ of $\dot{x}=f(x)$ such that $x_{z}(0)=z$ leaves $U_{\varepsilon}^{\circ}$
through of $S=U_{\varepsilon}\cap S_{\varepsilon}$. To see this, it's enough to show that $x_{z}$ can't remain in 
$U_{\varepsilon}$ in the future, since by $(ii)$, $(iii)$ and $(iv)$ we have that $x_{z}(t)\notin \partial U$ for all $t>0$.

Suppose, by contradiction, that this is false, then by using the compactness of $U_{\varepsilon}$ we have that $x_{z}$
is defined in $[0,\infty)$, and, since $W(0)=0$, we see by $(iii)$ and $(iv)$ that there is a $\delta>0$ such that
$\delta\leq |x_{z}(t)|\leq \varepsilon$, for all $t>0$.



Then, from $(i)$ we see that $h=V\circ x_{z}$ is a $\mathcal{C}^{1}$ function defined in $[0,\infty)$, and there is a $c>0$
such that $\dot{h}(t)<-c$.

This is a contradiction with the fact of $h$ is bounded in $[0,\infty)$ and shows that our assertion about $x_{z}$ is true.

Consider a sequence $x_{k}\in U_{\varepsilon}^{\circ}$ convergent to $0$, and denote by $\phi_{k}$
the solution of $\dot{x}=f(x)$ with $\phi_{k}(0)=x_{k}$. Let $y_{k}$ be the point of $S$ where $\phi_{k}$ leaves
$U_{\varepsilon}^{\circ}$, for the first time and $T_{k}>0$ the first positive number such that $\phi_{k}(T_{k})=y_{k}$.

Since $x_{k}\rightarrow 0$ which is an equilibrium of $\dot{x}=f(x)$ we have, by the continuous
dependence of the solutions of an o.d.e. with respect to the initial conditions, that
$T_{k}\rightarrow \infty$ when $k\rightarrow \infty$.

We can suppose, by the compactness of $S$, that $y_{k}$ converges to a point $\overline{y}\in S$. Let $\phi$ be the solution the consider differential equation with $\phi(0)=\overline{y}$.

Again by the continuous dependence of the solutions with respect to the initial conditions, it's easy to see that
there exists a $T>0$ such that $\phi(t)\in U_{\varepsilon}$ for $t\in [-T,0)$. 

Now we'll prove that $\phi(t)\in U_{\varepsilon}$
for all $t<0$ for which $\phi$ is defined.

If fact, suppose, by contradiction, that this is false then there exists a $\lambda<0$ such that
$\phi(\lambda)\notin U_{\varepsilon}$; since $T_{k}\rightarrow \infty$ there is a sufficiently large $k$ such that $T_{k}>2|\lambda|$.
Then, by the continuous dependence theorem, $\phi_{k}$ must leave $U_{\varepsilon}$ in a time $t_{k}$ such that
$0<t_{k}<T_{k}$. But this is a contradiction with the choice of $y_{k}$ and $T_{k}$.

Therefore, since $U_{\varepsilon}$ is compact, we have that $\phi$ is defined in $(-\infty,0]$ and the alpha limit set of this solution,
$\Lambda$, must be nonempty. Moreover, it's clear that $\Lambda\subset U_{\varepsilon}$.

By $(i)$, we have $\dot{V}\le 0$ in $U_{\varepsilon}$ therefore, by LaSalle's Invariance Principle, it follows that $\Lambda\subset\{q\in U_{\varepsilon}:\dot{V}(q)=0\}=\{0\}$ and we are done.
\edem

We make an imediate use of this to prove:

\begin{Proposition}\label{cor:hassymptotic}
Let $H=T+\Pi$ be a Hamiltonian function defined in $\Omega\times \mathbb{R}^{n}$, where the potential energy of class $\mathcal{C}^{2}$,
$\Pi:\Omega\rightarrow \mathbb{R}$, has a critical point at $0\in \Omega$. If there are $\varepsilon>0$ and a connected component $C$ of
$\Pi^{-1}((-\infty,0))$ adherent to the origin such that $\langle \nabla \Pi(q),q\rangle<0$ for
$q\in (\overline{C}\setminus \{0\})\cap B_{\varepsilon}$
then there is an asymptotic motion to $0$ in equations \eqref{eq1sec1}.
\end{Proposition}
\bdem\
It`s enough consider
\begin{eqnarray*}
U=\{(q,p)\in C\times \mathbb{R}^{n}|\quad \langle q,p\rangle >0,\quad |(q,p)|<\varepsilon \quad e\quad  H(q,p)<0 \}
\end{eqnarray*}
and functions $V,W: U\rightarrow \mathbb{R}$ defined by $V(q,p)=-\langle q,p\rangle$ and $W(q,p)=\langle q,p\rangle H(q,p)$.

Then, if we remember that $T$ is homogeneous of degree $2$ in relation to the variable $p$, a simple application of
Euler`s theorem for homogeneous functions and Hamiltonian equations shows that
\begin{eqnarray*}
\dot{V}(q,p)&=&-\Big[2T(q,p)-\Big\langle \frac{\partial T}{\partial q}(q,p),q \Big\rangle-\langle \nabla\Pi(q),q \rangle\Big],\\
\dot{W}(q,p)&=&\Big[2T(q,p)-\Big\langle \frac{\partial T}{\partial q}(q,p),q \Big\rangle-\langle \nabla\Pi(q),q \rangle\Big]H(q,p).
\end{eqnarray*}
Observe now that we can choose $\varepsilon$ small enough to yeld that
$R(q,p)=2T(q,p)-\langle \frac{\partial T}{\partial q}(q,p),q \rangle >0$ if $|q|<\varepsilon$ and $p\neq 0$. In fact, this is a
consequence of $R$ is quadratic with respect to $p$ and $R(0,p)=2T(0,p)$ is positive defined.

Then we see that conditions $(ii)$, $(iii)$ and $(iv)$ of proposition \ref{prop:karsovskii} are verified.

Now, if $(q,p)\in \overline{U}\setminus\{0\}$ we have $q\neq 0$, because if $p\neq 0$ then $H(0,p)>0$ and of course this
contradicts $(0,p)\in \overline{U}$.

This shows that $\dot{V}<0$ in $\overline{U}\setminus\{0\}$ and we can apply proposition \ref{prop:karsovskii}
in order to obtain our thesis.
\edem
Joining Theorem \ref{d5} and Proposition \ref{cor:hassymptotic} results immediately that if the potential energy of the Hamiltonian system is in the class $\mathcal{H}$, then there is an asymptotic trajectory to the origin.

\section{Proof of results of section \ref{sec2}}\label{sec.proofsomited}
\subsection{\hspace{2 cm} Proof of the Lemma \ref{b1} (and Corollary \ref{cor:compconexa})}
As $\Pi$ satisfies $(H_{3})$ there is $\varepsilon>0$ such that $\overline{B_{\varepsilon}}\subset\Omega$ and $(\overline{A}_{s}\setminus \{0\})\cap B_{\varepsilon}\subset C_{s}\cap B_{\varepsilon}$. Let $A_{*}$ be a connected component of $A_{s}$ such that $\overline{A}_{*}$ contains the immage of $\gamma$.

Let us fix some notation, take $q\in r_{\gamma}$, with $0<\Vert{q}\Vert<\varepsilon$ and, for $\rho>0$ denote $\Sigma_{\rho}(q)$ the disc of dimension $n-1$ centered in $q$ of radius $\rho$, orthogonal to $r_{\gamma}$. Since $\Vert{q}\Vert<\varepsilon$ we can choose $\rho$ small enough to guarantee  $\Sigma_{\rho}(q)\subset\overline{B_{\varepsilon}}$.

Consider the cone region $\Delta=\{z=\sigma w:\ 0<\sigma\le 1, w\in\Sigma_{\rho}(q)\}$. Then $\Delta\subset\overline{B_{\varepsilon}}$ and, as $r_{\gamma}$ is tangent to $\gamma$ in $0^{+}$, there is a $\delta>0$ such that $\gamma(t)\in\Delta$, if $0<t<\delta$.  Since for $t>0$, $\gamma(t)\in\overline{A}_{*}\setminus\{0\}$, it follows that exists $p\in\Delta\cap A_{*}$. Let $q_{\rho}$ the point of  $\Sigma_{\delta}(q)$ which is in the ray of origin $0$ which passes through $p$.

Clearly $j^{s}\Pi(p)<0$ and we claim that $j^{s}\Pi(x)<0$ for all $x$ of the segment $pq_{\rho}$. In fact, suppose by contradiction that this does not occur, then there is a point $y$ of segment $pq_{\rho}$ such that $j^{s}\Pi(x)<0$ for all $x$ between $p$ and $y$, and $j^{s}\Pi(y)=0$.

 Of course $y\in\overline{A}_{s}$, and as $p$ lies in the interior of segment $0q_{\rho}$, $y\not=0$ and follows from the previous considerations that $0\le\frac{d}{dt}j^{s}\Pi(ty)\vert_{t=1}=\langle\nabla (j^{s}\Pi)(y),y\rangle=R^{s}_{\Pi}(y)$.

Then $y\not=0$ and $y\in \overline{A}_{s}\setminus C_{s}$ against $(H_{3})$, which shows our claim.

Therefore the segment $pq_{q_{\rho}}$ is contained in $A_{*}$. This shows that, if $q\in r_{\gamma}$ with $0<\Vert{q}\Vert<\varepsilon$, then $q\in\overline{A_{*}}\setminus\{0\}$.

This shows the last affirmation of the lemma, moreover, by $(H_{3})$, one sees that $R^{s}_{\Pi}(q)<0$, for $q\in r_{\gamma}$ with $0<\Vert{q}\Vert<\varepsilon$. Since $r_{\gamma}$ is a ray and $j^{s}\Pi(0)=0$ this implies that $j^{s}\Pi(q)<0$ for these points. \edem

Corollary \ref{cor:compconexa} needn't an independent proof, the demonstration above does it.

\subsection{\hspace{2 cm} Proof of the proposition \ref{c2}}
By the curve selection lemma for semi-algebraic sets (see \cite{JM}) there is an algebraic curve $\gamma:[0,1]\longrightarrow\mathbb{R}^{n}$, with $\gamma(0)=0$ and $\gamma(t)\in A_{s}$ for $t>0$. Then, by lemma \ref{b1}, if $r$ is the tangent ray to $\gamma$ at $0+$, $j^{s}\Pi\vert_{r}$ has a strict maximum at origin. 

Therefore, by $(H_{1})$, we have $r\subset Z_{s-1}$ and $\Pi_{s}(q)<0$ for some $q\in r$, since $\Pi_{s}$ is an homogeneous polynomial of degree $s$, this shows that $r$ has the desired properties.

Let $\Delta_{r}$ be the connected component of $Z_{s-1}$ that contains $r$. We claim that $\Delta_{r}\subset A_{s}$. Indeed, note that if $x\in Z_{s-1}$ and $j^{s}\Pi(x)=0$, then $\Pi_{\ell}(x)=0$, for all $\ell\le s$, so $R^{s}_\Pi(x)=0$. Since $\Pi$ satisfies $(H_{3})$ this implies that $\partial A_{s}\setminus\{0\}\cap Z_{s-1}=\emptyset$. Our claim follows immediately from this observation, noting that
$\Delta_{r}$ is connected, $A_{s}$ is open, and $r\setminus\{0\}\subset A_{s}$.

To build the cone $K_{s}$, note that $Z_{s-1}$ is a cone with vertex $0$, so that $\Delta_{r}$ is also, then $\Delta_{r}$ is the union of half-lines of origin $0$ that are contained therein. Now take $\widetilde{r}$ a ray with origin $0$ contained in $\Delta_{r}$ and $x\in\widetilde{r}$, we have seen that $j^{s}\Pi(x)=\Pi_{s}(x)<0$, so we can choose $\delta>0$ such that, if $\Sigma_{\delta}(x)$ is the closed disc of dimension $n-1$ centred in $x$ and radius $\delta$, orthogonal to $\widetilde{r}$, then $\Pi_{s}(w)<0$ for all $w\in\Sigma_{\delta}(x)$. Moreover, because $\Delta_{r}$ is a connected component of cone $Z_{s-1}$, we can choose $\delta$ small enough in order to ensure that $\Sigma_{\delta}(x)\cap Z_{s-1}\subset\Delta_{r}$.

Consider the closed cone $K_{\widetilde{r}}=\{y\in\mathbb{R}^{n}:y=\lambda w,\ \lambda\ge 0,\ w\in\Sigma_{\delta}(x)$,
which has vertex at the origin and axis $\widetilde{r}$, since $\Pi_{s}$ is homogeneous, it's clear that $\Pi_{s}<0$ in $K_{\widetilde{r}}\setminus\{0\}$. Now define $K_{s}$ as
\begin{eqnarray*}
	K_{s}:=\bigcup_{\widetilde{r}\subset \Delta_{r}} K_{\widetilde{r}}
\end{eqnarray*}
Of course $K_{s}$ is a closed cone with vertex $0$ and $\Pi_{s}(q)<0$, if $q\in K_{s}\setminus\{0\}$, which shows item $(c)$.

As for every ray of origin $0$, $\widetilde{r}\subset \Delta_{r}$, we have $\widetilde{r}\subset K_{\widetilde{r}}^{\circ}\subset K_{s}^{\circ}$, $(a)$ follows.

In order to prove (b), observe that $K_{s}\cap Z_{s-1}=\bigcup_{\tilde{r}\subset \Delta_{r}} K_{\tilde{r}}\cap Z_{s-1}$ and, since
$K_{\tilde{r}}\cap Z_{s-1}\subset \Delta_{r}$, result $K_{s}\cap Z_{s-1}\subset \Delta_{r} \subset (K_{s})^{\circ}$.

As $\overline{Z_{s-1}}=Z_{s-1}\cup \{0\}$ is clear that, if $q\in \overline{Z_{s-1}}$ and $q\neq 0$, then $q\in Z_{s-1}$, which shows
that $\partial K_{s}\cap Z_{s-1}=\emptyset$. Therefore, if $q\in \partial K_{s}$ and $q\neq 0$, then $q\in Z_{s-1}^{c}$. 

So, since $	Z_{s-1}^{c}=\bigcup_{\ell=2}^{s-1}\big\{q\in \Omega : \Pi_{\ell}(q)\neq 0\big\}$, 
comes that, for $q\in \partial K_{s}\setminus \{0\}$, there is $\ell_{1}=\ell_{1}(q)\in \{2,...,s-1\}$ such that $\Pi_{\ell_{1}}(q)\neq 0$. 

Then, if $\ell(q):=\min\{\ell_{1}\in \{2,...,s-1\}:\Pi_{\ell_{1}}(q)\neq 0\}$, we have $(b)$.\edem

\subsection{\hspace{2 cm} Proof of the Proposition \ref{c3}}
Consider the cone $K_{s}$ built in the proposition \ref{c2}, clearly $\partial K_{s}$ is also a cone with vertex at the origin and, if $S$ is the sphere with center in the origin and radius $1$, $\partial K_{s}\cap S$ is compact.

For $x\in\partial K_{s}\setminus\{0\}$ consider the function given in item (b) of proposition \ref{c2}, $\ell(x)=\min\{\ell\in\mathbb{N}:\Pi_{\ell}(x)\not=0\}$, and take $F_{r}=\{x\in\partial K_{s}\setminus\{0\}:\ell(x)=r\}$, of course $\partial K_{s}\setminus\{0\}=\cup_{r=2}^{s-1}F_{r}$.

If $x\in \partial K_{s}\cap S$, item (b) of proposition \ref{c2} shows that $\Pi_{\ell}(x)>0$, for some $\ell\in\{2,\ldots,s-1\}$, so we can 
choose $\varepsilon_{x}>0$ such that $\Pi_{\ell}(q)>0$ for all $q\in (\partial K_{s}\cap S)\cap \overline{B_{\varepsilon_{x}}(x)}$. Then, put $B_{x}=(\partial K_{s}\cap S)\cap \overline{B_{\varepsilon_{x}}(x)}$ and  consider $m_{x}=\min_{q\in B_{x}}\{\Pi_{\ell}(q)\}$, obviously $m_{x}>0$.

If $\mathcal{K}(B_{x})$ is the cone with vertex in the origin generated by $B_{x}$, i.e. $\mathcal{K}(B_{x})=\{\lambda y:\lambda\ge 0,y\in B_{x}\}$, then for $q\in \mathcal{K}(B_{x})\setminus\{0\}$
\begin{eqnarray*}
	j^{\ell}\Pi(q)= j^{\ell-1}\Pi(q)+\Pi_{\ell}(q)\geq \Pi_{\ell}(q)=|q|^{\ell}\Pi_{\ell}\left(\frac{q}{|q|}\right)\geq m_{x}|q|^{\ell}.
\end{eqnarray*}
So, if $j^{s-1}_{\ell}\Pi=\Pi_{\ell+1}+...+\Pi_{s-1}$, comes 
\begin{eqnarray*}
	j^{s-1}\Pi(q)= j^{\ell}\Pi(q)+j^{s-1}_{\ell}\Pi(q)\geq m_{x}|q|^{\ell}+j^{s-1}_{\ell}\Pi(q),\ \forall q\in \mathcal{K}(B_{x}),\ q\neq 0,
\end{eqnarray*}
moreover, as $j^{\ell}(j^{s-1}_{\ell}\Pi)=0$, there is $\rho_{x}>0$ such that, if $0<|q|< \rho_{x}$, then $ |j^{s-1}_{\ell}\Pi(q)|\leq \frac{m_{x}}{2}|q|^{\ell} $

Thus, if $0<|q|< \rho_{x}$ and $q\in \mathcal{K}(B_{x})$, we have
\begin{eqnarray*}
	j^{s-1}\Pi(q)\geq \frac{m_{x}}{2}|q|^{\ell}
\end{eqnarray*}
The family $\mathcal{C}=\{B_{\varepsilon_{x}}(x):x\in \partial K_{s}\cap S\}$ is an open cover of $\partial K_{s}\cap S$.
By the compactness of $\partial K_{s}\cap S$, take $\{B_{\varepsilon_{x_{1}}}(x_{1}),B_{\varepsilon_{x_{2}}}(x_{2}),...,B_{\varepsilon_{x_{r}}}(x_{r})\}$
a finite subcover of $\mathcal{C}$, since $B_{x_{i}}=(\partial K_{s}\cap S)\cap \overline{B_{\varepsilon_{x_{i}}}(x_{i})}$, it's clear that
\begin{eqnarray*}
	\partial K_{s}\cap S \subset B_{x_{1}} \cup B_{x_{2}}\cup...\cup B_{x_{r}}.
\end{eqnarray*}
By the above procedure, for each $i\in \{1,...,r\}$ there are
$m_{x_{i}}>0$, $\rho_{x_{i}}>0$ and $\ell_{i}\in \{k,...,s-1\}$ such that, if $0<|q|<\rho_{x_{i}}$ and $q\in \mathcal{K}(B_{x_{i}})$ we have $j^{s-1}\Pi(q)\geq \frac{m_{x_{i}}}{2}|q|^{\ell_{i}}$.

Consider $
m:=\min\{m_{x_{1}},...,m_{x_{r}}\}>0,\ \rho:=\min\{\rho_{x_{1}},...,\rho_{x_{r}}\}>0$
and take $q\in \partial K_{s}$, $q\neq 0$. As $\partial K_{s}$ is a cone with vertex at
origin, $x_{q}=\frac{q}{|q|}\in \partial K_{s}\cap S$, hence
$x_{q}\in B_{\varepsilon_{x_{i_{q}}}}(x_{i_{q}})$ for some $i_{q}\in \{1,...,r\}$, consequently we have that 
$x_{q}\in B_{x_{i_{q}}}$, and so $q=|q|x_{q}\in \mathcal{K}(B_{x_{i_{q}}})$.

Thus, if $0<|q|<\rho$ we obtain $j^{s-1}\Pi(q)\geq \frac{m_{x_{i_{q}}}}{2}|q|^{\ell_{i_{q}}}\geq \frac{m}{2}|q|^{\ell_{i_{q}}}$,
therefore, if $0<\varepsilon^{\prime}<\min\{\rho,1\}$, it follows item (1) for $\varepsilon^{\prime}$, since
\begin{eqnarray*}
	j^{s-1}\Pi(q)\geq \frac{m}{2}|q|^{\ell_{i_{q}}}\geq \frac{m}{2}|q|^{s-1},\qquad \forall q\in (\partial K_{s}\setminus\{0\})\cap B_{\varepsilon^{\prime}} 
\end{eqnarray*}

To demonstrate (2), let $k$ be the degree of the first non null jet of $\Pi$, so $2\le k\le s-1$ and $j^{k}\Pi=\Pi_{k}$ is an homogeneous polynomial of degree $k$.

Note that $j^{s-1}\Pi(q)=\sum_{\ell=k}^{s-2}\Pi_{\ell}(q)+\Pi_{s-1}(q)$, so we have
\begin{center}
	$\Pi_{s-1}(q)\geq \frac{m}{2}|q|^{s-1}-\sum_{\ell=k}^{s-2}\Pi_{\ell}(q),\quad \forall q\in (\partial K_{s}\setminus\{0\})\cap
	B_{\varepsilon^{\prime}}$\ \ \ \ $(\dag)$ 
\end{center}
As $j^{\ell}\Pi(q)=\Pi_{\ell}(q)+\sum_{i=k}^{\ell-1}\Pi_{i}(q)$ and, by hypothesis $(H_{1})$,  $j^{\ell}\Pi(q)\geq 0$, if $k\le \ell \le s-1$, it follows
\begin{center}
	$\Pi_{\ell}(q)\geq -\sum_{i=k}^{\ell-1}\Pi_{i}(q),\quad \forall q\in B_{\varepsilon^{\prime}},\quad k+1\leq \ell \leq s-1.$ \ \ \ \ $(\ddag) $
\end{center}
The key point of this proof is to show, by induction, that if $2\le\nu\le s-k$ then for $0<t<\frac{k}{k+1}$, and $q\in (\partial K_{s}\setminus\{0\})\cap B_{\varepsilon^{\prime}}$,
\begin{center}
	$ P_{s-1}(tq)\ge\sum_{\ell=k}^{s-\nu}(\ell t^{\ell}-(s-\nu+1)t^{s-\nu+1})\Pi_{\ell}(q)+(s-1)t^{s-1}\frac{m}{2}\vert{q}\vert^{s-1}$\ \ \  $(*)$
\end{center}
Indeed, for $\nu=2$, by $(\dag)$,
$P_{s-1}(tq)=\sum_{\ell=k}^{s-2}\ell t^{\ell}\Pi_{\ell}(q)+(s-1)t^{s-1}\Pi_{s}(q)\ge\sum_{\ell=k}^{s-2}\ell t^{\ell}\Pi_{\ell}(q)+ (s-1)t^{s-1}(\frac{m}{2}\vert{q}\vert^{s-1}-\sum_{\ell=k}^{s-2}\Pi_{\ell}(q))$
and we have $(*)$.

Now suppose that $(*)$ worth for $2\le\nu<k+1$, and note that\newline $\sum_{\ell=k}^{s-\nu}(\ell t^{\ell}-(s-\nu+1)t^{s-\nu+1})\Pi_{\ell}(q)= \sum_{\ell=k}^{s-(\nu+1)}(\ell t^{\ell}-(s-\nu+1)t^{s-\nu+1})\Pi_{\ell}(q)+[(s-\nu)t^{s-\nu}-(s-\nu+1)t^{s-\nu+1}]\Pi_{s-\nu}(q)$.

As the sequence $(\frac{p}{p+1})$ is increasing and $0<t<\frac{k}{k+1}$, it follows that $(s-\nu)t^{s-\nu}-(s-\nu+1)t^{s-\nu+1}>0$, and the induction hypothesis together with $(\ddag)$ for $\Pi_{s-\nu}(q)$ show 
\begin{align*}
P_{s-1}(tq) & \ge
\sum_{\ell=k}^{s-(\nu+1)}(\ell t^{\ell} 
-(s-\nu+1)t^{s-\nu+1})\Pi_{\ell}(q)\cr
&-[(s-\nu)t^{s-\nu}-(s-\nu+1)t^{s-\nu+1}]\sum_{\ell=k}^{s-(\nu+1)}\Pi_{\ell}(q)+(s-1)t^{s-1}\frac{m}{2}\vert{q}\vert^{s-1} \cr
&=\sum_{\ell=k}^{s-(\nu+1)}(\ell t^{\ell}-(s-\nu+1)t^{s-\nu+1})\Pi_{\ell}(q)+(s-1)t^{s-1}\frac{m}{2}\vert{q}\vert^{s-1}
\end{align*}
This concludes the induction proof.

So, if $0<t<\frac{k}{k+1}$ and $q\in (\partial K_{s}\setminus\{0\})\cap B_{\varepsilon^{\prime}}$, then put $\nu=s-k$ in $(*)$ and conclude
$P_{s-1}(tq)\ge (k-(k+1)t)\Pi_{k}(tq)+(s-1)t^{s-1}\frac{m}{2}\vert{q}\vert^{s-1}$. 

Since $\Pi_{k}(tq)=j^{k}\Pi(tq)\geq 0$, if $0<t<\frac{k}{k+1}$, results from last inequality
\begin{eqnarray*}
	R^{s-1}(tq)\geq (s-1)\frac{m}{2}|tq|^{s-1},\quad \forall q\in (\partial K_{s}\setminus\{0\})\cap B_{\varepsilon^{\prime}},\quad 0<t<\frac{k}{k+1}
\end{eqnarray*}
Therefore, if $\varepsilon^{\prime\prime}=\frac{k}{k+1}\varepsilon^{\prime}<\varepsilon^{\prime}$ we obtain the stated in the item (2)
\begin{eqnarray*}
	R^{s-1}(q)\geq (s-1)\frac{m}{2}|q|^{s-1},\quad \forall q\in (\partial K_{s}\setminus\{0\})\cap B_{\varepsilon^{\prime\prime}}
\end{eqnarray*}
Then (1) and (2) are valid for $\varepsilon=\varepsilon^{\prime\prime}$.\edem

\subsection{\hspace{2 cm} Proof of the Proposition \ref{c5}}
Consider the closed cone $K_{s}$ with vertex at the origin built in proposition \ref{c2}, from part (c) of this proposition and homogeneity of $\Pi_{s}$, it follows that $\Pi_{s}(\frac{q}{|q|})=|q|^{-s}\Pi_{s}(q)<0$ for all $q\in K_{s}\setminus \{0\}$, as the sphere of radius $1$ is compact, there are strictly positive constants $M_{1}$ and $M_{2}$
such that
\begin{center}
	$-M_{1}\vert{q}\vert^{s}\le \Pi_{s}(q)\le
	-M_{2}\vert{q}\vert^{s},\quad \forall q \in (K_{s}\setminus \{0\}) $.
\end{center}
%
On the other hand, from Proposition \ref{c3},
there are $m>0$ and $0<\varepsilon^{\prime}<1$ such that
$ j^{s-1}\Pi(q)\geq \frac{m}{2}|q|^{s-1} $, for all 
$ q\in (\partial K_{s}\setminus\{0\})\cap B_{\varepsilon^{\prime}} $.

So, if $q\in (\partial K_{s}\setminus\{0\})\cap B_{\varepsilon^{\prime}}$, results
\begin{center}
	$j^{s}\Pi(q)=j^{s-1}\Pi(q)+\Pi_{s}(q)\geq \frac{m}{2}|q|^{s-1}-M_{1}|q|^{s}=|q|^{s-1}\Big(\frac{m}{2}-M_{1}|q|\Big)$.
\end{center}
This shows item $(i)$ of this proposition in $B_{\varepsilon_{1}}$, if $\varepsilon_{1}=\min\{\varepsilon^{\prime},\frac{m}{2M_{1}}\}$.

Again by Proposition \ref{c3}, we have  $0<\varepsilon^{\prime\prime}<1$ such that
\begin{center}
	${\displaystyle{R^{s-1}(q)\geq (s-1)\frac{m}{2}|q|^{s-1},\quad  \forall q\in (\partial K_{s}\setminus\{0\})\cap B_{\varepsilon^{\prime\prime}}}} $
\end{center}
then, for $q\in (\partial K_{s}\setminus\{0\})\cap B_{\varepsilon^{\prime\prime}}$, it follows from $R^{s}(q)=R^{s-1}(q)+s\Pi_{s}(q)$ that
\begin{center}
	${\displaystyle{R^{s}(q)\ge (s-1)\frac{m}{2}|q|^{s-1}-sM_{1}|q|^{s}= |q|^{s-1}\Big[(s-1)\frac{m}{2}-sM_{1}|q|\Big]}} $.
\end{center}

Then, if $\varepsilon_{2}=\min\{\varepsilon^{\prime\prime},\frac{(s-1)m}{2sM_{1}}\}$, we obtain item $(ii)$
in $B_{\varepsilon_{2}}$. 

To complete the proof take $\varepsilon=\min\{\varepsilon_{1},\varepsilon_{2}\}$.
\edem

\subsection{\hspace{2 cm} Proof of the Proposition \ref{c7}}
Let $\Delta_{r}$ be the connected component of $M_{s-1}$ which contains $r\setminus\{0\}$, follows from Corollary \ref{cor:compconexa} that $\Delta_{r}$ is the tangent cone to $\overline{A}_{*}\setminus\{0\}$ and $\Delta_{r}\subset A_{*}$. 

So, by part $(a)$ of Proposition \ref{c2}, $\Delta_{r}\subset (K_{s})^{\circ}$ and part $(b)$ of same proposition shows that $\partial K_{s}\setminus\{0\}$ is not tangent to $\overline{A}_{*}\setminus\{0\}$ at origin.


Moreover, item $(i)$ of proposition \ref{c5} implies that there is $\varepsilon_{1}>0$ such that $j^{s}\Pi(q)>0$ for all 
$q\in (\partial K_{s}\setminus\{0\})\cap B_{\varepsilon_{1}}$, so decreasing $\varepsilon_{1}$ if necessary, we obtain 	${\displaystyle{(\overline{A_{*}}\setminus \{0\})\cap B_{\varepsilon_{1}}\subset (K_{s})^{\circ}\cap B_{\varepsilon_{1}}}} $. 

As item $(c)$ of Proposition \ref{c2} states that $\Pi_{s}(q)<0$ for all $q\in K_{s}\setminus \{0\}$, it follows that, if $ q\in (\partial A_{*}\setminus \{0\})\cap B_{\varepsilon_{1}}$,
$  {\displaystyle{\Pi_{s}(q)<0}}$

By hypothesis $(H_{3})$, $(\overline{A}_{*}\setminus \{0\})\cap B_{\varepsilon_{1}}\subset
C_{*}\cap B_{\varepsilon_{1}}$, consequently  $R^{s}(q)<0$, for all 
$q\in (\partial A_{*}\setminus \{0\})\cap B_{\varepsilon_{1}}$.

Since $R^{s}(q)=\langle \nabla j^{s}\Pi(q),q\rangle=
(s-1)j^{s}\Pi(q)-\sum_{\ell=k}^{s-2}j^{\ell}\Pi(q)+\Pi_{s}(q)$, by using the hypothesis $(H_{1})$ we see that
\begin{center}
	$\displaystyle{R^{s}(q)\leq \Pi_{s}(q),\quad \forall q\in (\partial A_{*}\setminus \{0\})\cap B_{\varepsilon_{1}}}.\ \ \ \ \ \ \ (\dag) $
\end{center}


As $K_{s}$ is a positive cone with vertex in the origin, closed, it follows from the homogeneity of $\Pi_{s}$ and compacity of the sphere with center at origin and radius $1$ that there are $M_{1}>0$ and $M_{2}>0$ such that, for all $q \in K_{s}\setminus \{0\}$, 
$ \displaystyle{-M_{1}|q|^{s}\leq \Pi_{s}(q)\leq -M_{2}|q|^{s}} $.

As we had already proved that ${\displaystyle{(\overline{A_{*}}\setminus \{0\})\cap B_{\varepsilon_{1}}\subset (K_{s})^{\circ}\cap B_{\varepsilon_{1}}}} $, $(\dag)$ implies the desired inequality.

To complete the demonstration, use the item $(ii)$ of  Proposition \ref{c5} to find $\varepsilon_{2}>0$ such that $P_{s}(q)>0$ for all
$q\in (\partial K_{s}\setminus\{0\})\cap B_{\varepsilon_{2}}$.

Thus, if
$\varepsilon_{0}=\min\{\varepsilon_{1},\varepsilon_{2}\}$, we obtain $(\overline{C_{*}}\setminus \{0\})\cap B_{\varepsilon_{0}}\subset (K_{s})^{\circ}\cap B_{\varepsilon_{0}}$. \edem
\section{Final remarks}\label{sec:discuss}
Since the objective of this work is to find sufficient conditions to ensure that if $j^{s}\Pi$ is a potential of the
type Cetaev then the own potential energy $\Pi$ is a potential of the type Cetaev the hypothesis $(H_{2})$ and $(H_{3})$
are essential, without them visibly it is not possible to ensure this result.

Already the hypothesi $(H_{1})$ has a different character, it imposes a technical restriction that at least at first glance,
seems to be quite strong.

In this section we will show that some thecnical restriction need to be made on $j^{s}\Pi$, in addition to $(H_{2})$ and $(H_{3})$,
to ensure that $\Pi$ it is a potential of the type Cetaev, more precisely we will display an example that show that
$(H_{2})$ and $(H_{3})$ are not sufficient to ensure that $\Pi$ it is a potential of the type Cetaev, even in the case that 
$j^{s-1}\Pi(q)\geq 0$ in a neighborhood of the origin (which is a less restrictive restriction than that imposed by $(H_{1})$).

\begin{Exam}\label{d6}
	We consider the potential $\pi(x,y)=f(x,y)+x^{14}$, with
	\begin{center}
		$\displaystyle{f(x,y)=\frac{8}{3}y^{6}-3y^{4}x^{4}+\frac{9}{10}y^{2}x^{8}-\frac{1}{12}x^{12}-y^{12}}$.
	\end{center}
\end{Exam}
Then $j^{6}\pi$ is the first nonzero jet of $\pi$ and it is positive semi-definite. Also note that
$j^{12}\pi=f$ shows that $\pi$ has no minimum at the origin, because $f(x,0)=-\frac{x^{12}}{12}$.


So, $\pi$ obeys hypothesis $(H_{2})$ for $s=12$, furthermore, 
\begin{eqnarray*}
	j^{11}\pi(x,y)=\frac{8}{3}y^{6}-3y^{4}x^{4}+\frac{9}{10}y^{2}x^{8}=y^{2}(\frac{8}{3}y^{4}-3y^{2}x^{4}+\frac{9}{10}x^{8})
\end{eqnarray*}
is positive semi-definite and vanishes only on the axis $\{y=0\}$. The $j^{12}\pi$ is the first jet of $\pi$ to shows that this potential hasn't a minimum at origin.

However,
\begin{eqnarray*}
	j^{8}\pi(x,y)=\frac{8}{3}y^{6}-3y^{4}x^{4}=
	\frac{8}{3}y^{4}\Big(y-\frac{3}{2\sqrt{2}}x^{2}\Big)\Big(y+\frac{3}{2\sqrt{2}}x^{2}\Big)
\end{eqnarray*}
has a saddle at the origin.

This shows that $\pi$ does not obey $(H_{1})$, moreover it points clear an interesting fact, if the first nonzero jet of $\pi$ is positive semi-definite, $j^{s}\pi$ is the first jet that shows that
$\pi$ has no minimum at the origin and $j^{s-1}\pi$ is positive semi-definite, we can't conclude that $j^{\ell}\pi \geq 0$ if
$\ell \leq s-1$.

In order to prove that $\pi$ satisfies $(H_{3})$ for $s=12$ consider the radial derivative of $f$,
$R_{f}(x,y)=y^{6}-24y^{4}x^{4}+9y^{2}x^{8}-x^{12}-12y^{12}$ and note that, if $Q(x,y)=(y-x^{2})(y+x^{2})(2y-x^{2})^{2}(2y+x^{2})^{2}$, we have $R_{f}(x,y)=Q(x,y)-y^{12}$, so $R_{f}(x,y)\le Q(x,y)$.

As $Q(x,y)\le 0$ in $\Delta=\{(x,y)\in \mathbb{R}^{2}| -x^{2}\leq y \leq x^{2}\}$, we have that the connected component of $R_{f}^{-1}(]-\infty,0[)$ that contains the axis $y=0$ contains also $\Delta$.

Since $f(x,\pm x^{2})=f(x,\pm x^{2})=\left(\frac{29}{60}-x^{12}\right)x^{12}$, it follows that, for $\varepsilon=\sqrt[12]{\frac{29}{60}}$, in $B_{\varepsilon}$ the connected component of $f^{-1}(]-\infty,0[)$ which contains the axis $y=0$ is contained in $\Delta$. Therefore $\pi$ satisfies $(H_{3})$ for $s=12$.

Then, for $s=12$ we have that $\pi$ verifies $(H_{2}),\ (H_{3})$, and $j^{s-1}\pi$ has a nonstrict minimum at origin.  However, we claim that $\pi$ isn't a Cetaev potential.

To see this, note that $j^{12}\pi=f$ and $f(x,\lambda x^{2})= h(\lambda)x^{12}-(\lambda x^{2})^{12}\leq h(\lambda)x^{12}$, where $h(\lambda)=\frac{8}{3}\lambda^{6}-3\lambda^{4}+\frac{9}{10}\lambda^{2}-\frac{1}{12}$.

A direct calculus shows that $h(\lambda)<0$ in $[-\frac{3}{4},\frac{3}{4}]$, so there is a $\rho>0$ such that $f(x,\lambda x^{2})\le -\rho^{2} x^{12}$, if $\vert{\lambda}\vert\le \frac{3}{4}$.

Therefore, if $\vert{\lambda}\vert\le \frac{3}{4}$, $\pi(x,\lambda x^{2})=f(x,\lambda x^{2})+x^{14}\le (-\rho^{2} +x^{2})x^{12}$, and this shows that $\pi(x,y)<0$, if $(x,y)\in\{(x,y)\in \mathbb{R}^{2}: -\frac{3}{4} x^{2}\le y \le \frac{3}{4} x^{2}\}\cap B_{\rho}\setminus\{0\}$.

Since the radial derivative of $\pi$ is $R_{\pi}(x,y)=R_{f}(x,y)+14x^{14}$, it's clear that $R_{\pi}(x,\frac{x^{2}}{2})=-\frac{2}{2^{12}}x^{24}+14x^{14}$ is positive, if $\vert{x}\vert$ is small enough, so $\pi$ isn't a Cetaev potential, as we assert.


\end{document}